\documentclass[12pt]{article}
\usepackage{amsfonts,amsmath,amssymb,amsbsy,
amsthm,latexsym,amsopn,amstext,amsxtra,euscript,amscd,mathrsfs}

\newtheorem{theorem}{Theorem}
\newtheorem{lemma}{Lemma}

\newtheorem{proposition}{Proposition}
\theoremstyle{definition}
\newtheorem{definition}{Definition}

\def\Fp{\mathbb{F}_p}
\def\Fq{\mathbb{F}_q}

\begin{document}

\title{Average estimate for additive energy in prime field.}
\author{Glibichuk Alexey\thanks{Technion, Israel Institute of Technology, Haifa, Israel.
\endgraf E-mail:glibichu@tx.technion.ac.il.}}
\date{}

\maketitle

\begin{abstract}  Assume that $A\subseteq \Fp, B\subseteq \Fp^{*}$,
$\frac{1}{4}\leqslant\frac{|B|}{|A|},$ $|A|=p^{\alpha},
|B|=p^{\beta}$. We will prove that for $p\geqslant p_0(\beta)$ one
has
$$\sum_{b\in B}E_{+}(A, bA)\leqslant 15 p^{-\frac{\min\{\beta, 1-\alpha\}}{308}}|A|^3|B|.$$
Here $E_{+}(A, bA)$ is an additive energy between subset $A$ and
it's  multiplicative shift $bA$. This improves previously known
estimates of this type.
\end{abstract}

\section{Introduction.}

Let X be a non-empty set endowed with a binary operation $*:X\times
X\rightarrow X$. Then one can define the operation * on pairs of
subsets $A,B\subset X$ by the formula $A*B=\{a*b : a\in A, b\in
B\}.$ In particular, if $A$ and $B$ are subsets of a ring, we have
two such operations: addition $A+B:=\{a+b:a\in A, b\in B\}$ and
multiplication $AB=A\times B:=\{ab:a\in A, b\in B\}.$ For given
element $b$ we define operation $b*A={b}\times A$.  The sign $*$ may
be omitted when there is no danger of confusion. We write $|A|$ for
the cardinality of $A$. We take the ring to be the field $\Fp$ of
$p$ elements, where $p$ is an arbitrary prime. All sets are assumed
to be subsets of $\Fp$. Given any set $Y\subset\Fp$, we write
$Y^{*}:=Y\setminus\{0\}$ for the set of invertible elements of $Y$.
We shall always assume that $p$ is a prime. Given any real number y,
we write $[y]$ for its integer part (the largest integer not
exceeding y), and denote the fractional part of y by $\{y\}$. We
also define the operation $h+A=\{h\}+A$ which adds an arbitrary
element $h\in\Fp$ to the set $A$.

\begin{definition} For subsets $A, B\subset\Fp $ we denote
$$E_{+}(A,B)=|\{(a_1,a_2,b_1,b_2)\in A\times A\times B\times B: a_1-a_2=b_1-b_2\}|,$$
$$E_{\times}(A,B)=|\{(a_1,a_2,b_1,b_2)\in A\times A\times B\times B: a_1a_2=b_1b_2\}|.$$
Numbers $E_{+}(A,B)$ and $E_{\times}(A,B)$ are said to be an
\textbf{additive energy} and a \textbf{multiplicative energy} of
sets $A$ and $B$ respectively.
\end{definition}

In the paper \cite{B} J. Bourgain proved the following result.

\begin{theorem}\label{Bourgain}
Assume $A\subset\Fp, B\subset\Fp$ and $|A|=p^{\alpha},
|B|=p^{\beta}$ with $\alpha\geqslant\beta$. Then
$$\sum_{b\in B}E_{+}(A,bA)<C_1p^{c_2\gamma}|A|^3|B|$$
where $\gamma=\min(\beta,1-\alpha)$ and $C_1,c_2$ are absolute
constants (independent on $\alpha,\beta$).
\end{theorem}
In the same paper J. Bourgain  deduces from Theorem \ref{Bourgain}
sum-product estimate for two different subsets. Further, J. Bourgain
and author \cite{BoGl} of this paper extended Theorem \ref{Bourgain}
to the case of an arbitrary finite field. More precisely, we proved
the following result.
\begin{theorem}\label{BouGlib}Take arbitrary subsets $A,B$ of a finite field
$\Fq$ with $q=p^r$ elements, such that $|A|=q^{\alpha},
|B|=q^{\beta}, \alpha\geqslant\beta$ and an arbitrary
$0<\eta\leqslant 1$. Suppose further that for every nontrivial
subfield $S\subset\Fq$ and every element $d\in Fq$ the set $B$
satisfies the restriction
$$|B\cap dS|\leqslant 4|B|^{1-\eta}.$$
Then
$$\sum_{b\in B}E_{+}(A,bA)\leqslant 13q^{-\frac{\gamma}{10430}}|A|^{3}|B|$$
where
$\gamma=\min\left(\beta,\frac{5215}{4}\beta\eta,1-\alpha\right)$.
\end{theorem}
In this paper we also deduced from the Theorem \ref{BouGlib} a new
character sum estimate over a small multiplicative subgroup. J.
Bourgain, S. J. Dilworth, K. Ford, S. Konyagin and D. Kutzarova
\cite{BDFKK} applied Theorem \ref{BouGlib} to one of the problems of
sparse signal recovery and several others branches of coding theory.
Also, M. Rudnev and H. Helfgott \cite{RH} used method, proposed in
the proof of the Theorem \ref{Bourgain} to obtain an new explicit
point-line incidence result in $\Fp$. These examples demonstrate
that estimates like Theorems \ref{Bourgain} and \ref{BouGlib} have
wide range of applications.

In the current paper a slightly modified version of the method from
paper \cite{RH} will be used to obtain an improvement of the Theorem
\ref{BouGlib} in the case of prime field $\Fp$. We will establish
the following theorem.
\begin{theorem}\label{main1}  Assume that $A\subseteq \Fp, B\subseteq \Fp^{*}$,
$\frac{1}{4}\leqslant\frac{|B|}{|A|},$ $|A|=p^{\alpha},
|B|=p^{\beta}$. Then for $p\geqslant p_0(\beta)$
$$\sum_{b\in B}E_{+}(A, bA)\leqslant 15 p^{-\frac{\min\{\beta, 1-\alpha\}}{308}}|A|^3|B|.$$
\end{theorem}

Ideas of M. Rudnev and H. Helfgott in context of this problem
working only when $|B|\geqslant K|A|$ for some absolute constant
$K$. Case when $|A|$ is small comparatively to $|B|$ was analyzed by
another method. This method is elementary in some extent and gives
the following estimate.
\begin{theorem}\label{mainLargeB} Assume that $A\subseteq \Fp, B\subseteq \Fp^{*}$,
$|A|=p^{\alpha}, |B|=p^{\beta}$. Then for $p\geqslant p_0(\alpha,
\beta)$ we have
$$\sum_{b\in B}E_{+}(A, bA)\leqslant C p^{-\frac{\min\{\beta, 1-\alpha\}}{2240}}|A|^3|B|,$$
where $C>0$ is an absolute constant.
\end{theorem}
As we see, Theorem \ref{mainLargeB} gives worse estimate than
Theorem \ref{main1}, but it still better than one delivered by the
Theorem \ref{BouGlib}.

In section \ref{Sec1} we stating preliminary results which will be
used in proofs of Theorems \ref{main1} and \ref{mainLargeB}. Theorem
\ref{main1} is proved in the Section \ref{sec2}, Theorem
\ref{mainLargeB} is proved in the Section \ref{sec3}.

\textbf{Acknowledgements.} The author thank professor S. Konyagin
and M. Rudnev for useful discussions helped me to improve the final
result.

\section{Preliminary results.}\label{Sec1}

All the subsets in the Lemmas below are assumed to be non-empty. The
first two lemmas is due to Ruzsa \cite{R1, R2}. It holds for subsets
of any abelian group, but here we state them only for the subsets of
$\Fp$.

\begin{lemma}\label{RuszaTriangle} For any subsets $X$, $Y$, $Z$ of $\Fp$ we have
$$|X-Z|\leqslant\frac{|X-Y||Y-Z|}{|Y|}.$$
\end{lemma}

\begin{lemma}\label{RuszaGeneral} Let $Y,X_1,X_2,\ldots,X_k$ be sets
of $\Fp$. Then
$$|X_1+X_2+\ldots+X_k|\leqslant\frac{\prod_{i=1}^{k}|Y+X_{i}|}{|Y|^{k-1}}.$$
\end{lemma}

\begin{definition} For any nonempty subsets $A\subset\Fp, B\subset\Fp, G\subset
A\times B$, we define their \emph{partial sum}
$$|A\substack{+\\G}B|=\{a+b:(a,b)\in G\}.$$
\end{definition}

Let us recall the modification of Balog-Szemeredi-Gowers result (see
the paper of J. Bourgain and M. Garaev \cite{BG}, Lemma 2.3).

\begin{proposition}\label{BSG} Let $A$ and $B$ be subsets of $\Fp$
and $G\subset A\times B$ be such that $|G|\geqslant\frac{|A||B|}{K}$
for some $K>0$. Then there exist subsets $A^{'}\subset A,
B^{'}\subset B$ and a number $Q$, with
$$|A^{'}|\geqslant\frac{|A|}{4\sqrt{2}K},\phantom{111}\frac{|A|}{8\sqrt{2}K^2\ln(e|A|)}
\leqslant Q\leqslant 2|A^{'}|,\phantom{111}|B^{'}|\geqslant
\frac{|A||B|}{8\sqrt{2}QK^2\ln(e|A|)}$$ such that
$$|A\substack{+\\G}B|^3\geqslant |A^{'}+B^{'}|\frac{Q|B|}{256K^3\ln(e|A|)}.$$
\end{proposition}

We shall use the following result from the book of T. Tao and V. Vu
\cite{TV} (Lemma 2.30, p. 80).

\begin{lemma}\label{TaoVuLemma} If $E_{+}(A,B)>\frac{1}{K}|A|^{\frac{3}{2}}|B|^{\frac{3}{2}}, K\geqslant 1$,
then there is $G\subset A\times B$ satisfying
$$|G|>\frac{1}{2K}|A||B|\phantom{1}\textmd{and}\phantom{1}|A\substack{+\\G}B|<2K|A|^{\frac{1}{2}}|B|^{\frac{1}{2}}.$$
\end{lemma}

This lemma represents a known technical approach for estimating
sum-product sets, see, for example \cite{BKT}, \cite{BK}.

\begin{lemma}\label{XYGlemma} For any given subsets $X, Y\subseteq\Fp, G\subset\Fp^{*}$
there is an element $\xi\in G$ with
$$
|X+\xi Y|\geqslant\frac{|X||Y||G|}{|X||Y|+|G|}.
$$
Moreover, the following inequality holds
$$|X+\xi Y|>\frac{|X|^2|Y|^2}{E_{+}(X,\xi Y)}.$$
\end{lemma}

\textbf{Proof.} Let us take an arbitrary element $\xi\in G$ and
$s\in\Fp$ and denote
$$f_{\xi}^{+}(s):=|\{(x,y)\in X\times Y:x+y\xi =s\}|.$$
It is obvious that
$$\sum_{s\in
\Fp}(f_{\xi}^{+}(s))^{2}=|\{(x_1,y_1,x_2,y_2)\in X\times X\times
Y\times Y: x_1+y_1\xi = x_2+y_2\xi \}|$$
$$=|X||Y|+|\{(x_1,y_1,x_2,y_2)\in X\times X\times
Y\times Y: x_1 \neq x_2,x_1+y_1\xi=x_2+y_2\xi \}|
$$
and
\begin{equation}\label{XYGtmplab1}
\sum_{s\in\Fp}f_{\xi}^{+}(s)=|X||Y|.
\end{equation}

Let us observe that for every $x_1,x_2\in X,y_1,y_2\in Y$ such that
$x_1\neq x_2$, there is at most one $\eta\in G$ satisfying the
equality $x_1+y_1\eta=x_2+y_2\eta$. Therefore,
$$\sum_{\xi\in G}\sum_{s\in\Fp}(f_{\xi}^{+}(s))^{2}\leqslant
|X||Y||G|+|X|^{2}|Y|^{2}.$$ From the last inequality it directly
follows that there is an element $\xi\in G$ such that
\begin{equation}\label{XYGtmplab2}
\sum_{s\in\Fp}(f_{\xi}^{+}(s))^{2}\leqslant
|X||Y|+\frac{|X|^{2}|Y|^{2}}{|G|}.
\end{equation}
According to Cauchy-Schwartz,
\begin{equation}\label{XYGtmplab3}
\left(\sum_{s\in\Fp}f_{\xi}^{+}(s)\right)^{2}\leqslant |X+\xi
Y|\sum_{s\in\Fp}(f_{\xi}^{+}(s))^{2}.
\end{equation}
Observing that
$$\sum_{s\in
\Fp^{*}}(f_{\xi}^{+}(s))^{2}=E_{+}(X,\xi Y)$$ one can yield the
second assertion of Lemma \ref{XYGlemma}.

Combining inequalities (\ref{XYGtmplab1}), (\ref{XYGtmplab2}) and
(\ref{XYGtmplab3}) we see that
$$|X+\xi Y|\geqslant
\frac{|X|^2|Y|^2}{|X||Y|+\frac{|X|^2|Y|^2}{|G|}}=\frac{|X||Y||G|}{
|X||Y|+|G|}.
$$ Lemma \ref{XYGlemma} now follows. $\blacksquare$

\begin{definition} For any given subsets $X, Y\subset\Fp, |Y|>1$ we
denote
$$Q[X,Y]=\frac{X-X}{(Y-Y)\setminus\{0\}}:=\left\{\frac{x_1-x_2}{y_1-y_2}:
x_1,x_2\in X, y_1,y_2\in Y, y_1\neq y_2\right\}.$$ If $X=Y$ then
$Q[X,X]=Q[X]$.
\end{definition}

Lemma \ref{QXYproplemma} is a simple extension of Lemma 2.50 from
the book by T. Tao and V. Vu \cite{TV}.

\begin{lemma}\label{QXYproplemma} Consider two arbitrary subsets $X, Y\subset\Fp,
|Y|>1$. The given element $\xi\in\Fp$ is contained in $Q[X,Y]$ if
and only if $|X+\xi* Y|<|X||Y|$.
\end{lemma}
\textbf{Proof.} Let us consider a mapping $F:X\times Y$ to $X+\xi*
Y$ defined by the identity $F(x,y)=x+\xi y$. $F$ can be
non-injective only when $|X+\xi* Y|<|X||Y|$. On the other side, the
non-injectivity of $F$ means that there are elements $x_1,x_2\in X$,
$y_1, y_2\in Y$ such that $(x_1,y_1)\ne (x_2,y_2)$ and
$F(x_1,y_1)=F(x_2,y_2)$. It is obvious that $y_1\ne y_2$ since
otherwise $x_1=x_2$ and we have achieved a contradiction with
condition $(x_1,y_1)\ne (x_2,y_2)$. Hence,
$\xi=(x_1-x_2)/(y_2-y_1)\in Q[X,Y]$. Lemma \ref{QXYproplemma} now
follows. $\blacksquare$
\bigskip

We need the following Lemma due to C.-Y. Shen \cite{Shen}.

\begin{lemma}\label{ChenCovering} Let $X_1$ and $X_2$ be two sets.
Then for any $\varepsilon\in (0,1)$ there exist at most
$\frac{\ln{\frac{1}{\varepsilon}}}{|X_2|}\min \left\{|X_1+X_2|,
|X_1-X_2|\right\}$ additive translates of $X_2$ whose union contains
not less than $(1-\varepsilon)|X_1|$ elements of $X_1$.
\end{lemma}
\textbf{Proof.} For simplicity, we assume that $|X_1+X_2|\leqslant
|X_1-X_2|$. The case when $|X_1+X_2|>|X_1-X_2|$ can be considered
similarly. Using Lemma \ref{XYGlemma} we deduce
$$|\{(x,y,x_1,y_1)\in X_1\times X_2\times X_1\times X_2: x+y=x_1+y_1\}|
\geqslant\frac{|X_1|^2|X_2|^2}{|X_1+X_2|}.$$ Now we can fix two
elements $x_{*}^{1}\in X_1, y_{*}^{1}\in X_2$ for which the equation
$x_{*}^{1}+y=x+y_{*}^{1}, x\in X_1, y\in X_2$ has at least
$\frac{|X_1||X_2|}{|X_1+X_2|}$ solutions and, therefore,
$|(x_{*}^{1}+X_2)\cap
(y_{*}^{1}+X_1)|\geqslant\frac{|X_1||X_2|}{|X_1+X_2|}.$ Denoting
$K=\frac{|X_1+X_2|}{|X_2|}$ we can observe that
\begin{equation}\label{temp1}
|X_1\cap (x_{*}^{1}-y_{*}^{1}+X_2)|\geqslant\frac{|X_1|}{K}.
\end{equation}
Obviously, from (\ref{temp1}) it is  follows that
$$|X_{1}^{1}|:=|X_1\setminus
(x_{*}^{1}-y_{*}^{1}+X_2)|\leqslant\left(1-\frac{1}{K}\right)|X_1|.$$
We can repeat previous arguments for sets $X_1^{1}$ and $X_2$ and
find elements $x_{*}^{2}\in X_1^{1}$ and $y_{*}^{2}\in X_2$ such
that
$$|X_1^{1}\cap (x_{*}^{2}-y_{*}^{2}+X_2)|\geqslant\frac{|X_1^{1}|}{K}$$
$$|X_1^{2}|:=|X_1^{1}\setminus (x_{*}^{2}-y_{*}^{2}+X_2)|\leqslant
\left(1-\frac{1}{K}\right)|X_1^{1}|\leqslant
\left(1-\frac{1}{K}\right)^{2}|X_1|.$$ On $i$-th iteration we
finding elements $x_{*}^{i}\in X_1^{i-1}$ and $y_{*}^{i}\in X_2$
with
$$|X_1^{i-1}\cap (x_{*}^{i}-y_{*}^{i}+X_2)|\geqslant\frac{|X_1^{i-1}|}{K}$$
$$|X_1^{i}|:=|X_1^{i-1}\setminus (x_{*}^{i}-y_{*}^{i}+X_2)|\leqslant
\left(1-\frac{1}{K}\right)|X_1^{i-1}|\leqslant
\left(1-\frac{1}{K}\right)^{i}|X_1|.$$ We stop when
$|X_1^{n}|<\varepsilon |X_1|$ for some $n$. It is easy to see that
we will make not more than $\ln\left(\frac{1}{\varepsilon}\right)K$
steps. The last observation finishes the proof of the Lemma
\ref{ChenCovering}. $\blacksquare$

We also need the following sum-product estimate of M. Z. Garaev
\cite[Theorem 3.1]{Gar}.

\begin{theorem}\label{garaev} Let $A, B\subset\Fp$ be an arbitrary
subsets. Then
$$|A-A|^2\cdot\frac{|A|^2|B|^2}{E_{\times}(A,B)}\geqslant C
|A|^3L^{\frac{1}{9}}(\log_2L)^{-1},$$ where
$L=\min\left\{|B|,\frac{p}{|A|}\right\}$ and $C>0$ is an absolute
constant.
\end{theorem}

\section{Proof of the Theorem \ref{main1}.}\label{sec2}

Let $A,B\subseteq\Fp$ be as in Theorem \ref{main1} and $\delta>0$,
$C>1$ (to be specified). Assume
$$\sum_{b\in B}E_{+}(A, bA)>C|B|^{1-\delta}|A|^3.$$
Hence there is a subset $B_1\subseteq B$ such that
$$|B_1|>\frac{C}{2}|B|^{1-\delta}$$
and
\begin{equation}\label{tmplab3}
E_{+}(A,
bA)>\frac{C}{2}|B|^{-\delta}|A|^3\phantom{1}\textmd{for}\phantom{1}
b\in B_1.
\end{equation}

Fix $b\in B_1.$ By the application of Lemma \ref{TaoVuLemma} to
(\ref{tmplab3}), one can deduce that there is $G^{(b)}\subset
A\times bA, |G^{(b)}|>\frac{C}{4}|B|^{-\delta}|A|^2$ such that
$$|A\substack{+\\G^{(b)}}bA|<\frac{4}{C}|B|^{\delta}|A|.$$
Now, by Proposition \ref{BSG}, there are $Q_{(b)}, A_1^{(b)},
A_2^{(b)}\subset A$ such that
\begin{equation}\label{tmplabA1}
|A_{1}^{(b)}|>\frac{C}{2^{4}\sqrt{2}}|B|^{-\delta}|A|,
\end{equation}
\begin{equation}\label{tmplabQb}
\frac{C^2}{2^{7}\sqrt{2}\ln(e|A|)}|A||B|^{-2\delta}\leqslant Q_{(b)}
\leqslant 2|A_1^{(b)}|,
\end{equation}
\begin{equation}\label{tmplabA2}
|A_{2}^{(b)}|>\frac{C^2}{2^{7}\sqrt{2}Q_{(b)}\ln(e|A|)}|B|^{-2\delta}|A|^2,
\end{equation}
\begin{equation}\label{tmplab5}
|A_1^{(b)}+bA_2^{(b)}|<\frac{2^{20}}{C^6Q_{(b)}}\ln(e|A|)|B|^{6\delta}|A|^2.
\end{equation}
Write
\begin{multline*}
\frac{C^3}{2^{12}\ln(e|A|)}|B_1||B|^{-3\delta}|A|^2<\sum_{b\in
B_1}|A_1^{(b)}\times A_2^{(b)}|\\ \leqslant |A|\left[\sum_{b,
b^{'}\in B_1}\left|\left(A_1^{(b)}\cap A_1^{(b^{'})}\right)\times
\left(A_2^{(b)}\cap A_2^{(b^{'})}\right)\right|\right]^{\frac{1}{2}}
\end{multline*}
by Cauchy-Schwartz. Hence
$$\frac{C^6}{2^{24}\ln^2(e|A|)}|B_1|^2|B|^{-6\delta}|A|^2<
\sum_{b, b^{'}\in B_1}\left|\left(A_1^{(b)}\cap
A_1^{(b^{'})}\right)\times \left(A_2^{(b)}\cap
A_2^{(b^{'})}\right)\right|$$ and there is some $b_0\in B_1,
B_2\subset B_1$ such that
\begin{equation}\label{tmplabnew1}
|B_2|>\frac{C^7}{2^{26}\ln^2(e|A|)}|B|^{1-7\delta}
\end{equation}
\begin{equation}\label{tmplab7}
|A_1^{(b)}\cap A_1^{(b_0)}|, |A_2^{(b)}\cap A_2^{(b_0)}|>
\frac{C^6}{2^{25}\ln^{2}(e|A|)}|B|^{-6\delta}|A|\phantom{1}\textmd{for}\phantom{1}b\in
B_2.
\end{equation}
Let us estimate from (\ref{tmplabA1}), (\ref{tmplabA2}),
(\ref{tmplab5}), (\ref{tmplab7}) and Lemma \ref{RuszaTriangle}
\begin{multline}\label{multline1}
|b_0A_1^{(b_0)}+bA_1^{(b_0)}|\leqslant\frac{
|A_1^{(b_0)}+bA_2^{(b_0)}||A_1^{(b_0)}+b_0A_2^{(b_0)}|}{|A_2^{(b_0)}|}\leqslant\\
\leqslant\frac{2^{27}\sqrt{2}\ln^2(e|A|)}{C^{8}}|B|^{8\delta}|A_1^{(b_0)}+bA_2^{(b_0)}|
\end{multline}
\begin{multline}\label{multline2}
|A_1^{(b_0)}+bA_2^{(b_0)}|\leqslant\frac{
|A_1^{(b_0)}+bA_2^{(b)}||A_2^{(b_0)}+A_2^{(b_0)}|}{|A_2^{(b)}\cap
A_2^{(b_0)}|}\leqslant\\ \leqslant\frac{
|A_1^{(b_0)}+bA_2^{(b)}||A_1^{(b_0)}+b_0A_2^{(b_0)}|^2}{|A_2^{(b)}\cap
A_2^{(b_0)}||A_1^{(b_0)}|}\leqslant\\
\leqslant\frac{2^{69}\sqrt{2}\ln^{4}(e|A|)}{C^{19}Q^{2}_{(b_0)}}
|A|^2|B|^{19\delta}|A_1^{(b_0)}+bA_2^{(b)}|
\end{multline}
\begin{multline}\label{multline3}
|A_1^{(b_0)}+bA_2^{(b)}|\leqslant
\frac{|A_1^{(b)}+bA_2^{(b)}||A_1^{(b_0)}+A_1^{(b_0)}|}{|A_1^{(b_0)}\cap
A_1^{(b)}|}\leqslant\\ \leqslant
\frac{|A_1^{(b)}+bA_2^{(b)}||A_1^{(b_0)}+b_0A_2^{(b_0)}|^2}{|A_1^{(b_0)}\cap
A_1^{(b)}||A_2^{(b_0)}|}\leqslant \\ \leqslant
\frac{2^{92}\sqrt{2}\ln^6(e|A|)}{C^{26}Q_{(b)}Q_{(b_0)}}|B|^{26\delta}|A|^3.
\end{multline}
Hence, by (\ref{multline1}), (\ref{multline2}) and (\ref{multline3})
$$
|b_0A_1^{(b_0)}+bA_1^{(b_0)}|\leqslant
\frac{2^{189}\sqrt{2}\ln^{12}(e|A|)}{C^{53}Q_{(b_0)}^{3}Q_{(b)}}|B|^{53\delta}|A|^5.
$$
Using (\ref{tmplabQb}) finally we obtain
$$|b_0A_1^{(b_0)}+bA_1^{(b_0)}|\leqslant
\frac{2^{219}\sqrt{2}\ln^{16}(e|A|)}{C^{61}}|B|^{61\delta}|A|.$$

Now we redefine $A_1^{(b_0)}$ by $A^{'}$ and $\frac{B_2}{b_0}$ by
$B^{'}$ one can deduce the following properties (for
$\delta<\frac{1}{440}$):
\begin{equation}\label{AplbA}
|A^{'}+bA^{'}|<\frac{2^{219}\sqrt{2}\ln^{16}
(e|A|)}{C^{61}}|B|^{61\delta}|A|\phantom{1}\textmd{for}\phantom{1}\textmd{all}\phantom{1}b\in
B^{'}
\end{equation}
\begin{equation}\label{Bprime}
|B^{'}|>\frac{C^7}{2^{26}\ln^2(e|A|)}|B|^{1-7\delta}
\end{equation}
\begin{equation}\label{Aprime}
|A^{'}|>\frac{C}{2^{4}\sqrt{2}}|B|^{-\delta}|A|.
\end{equation}
Our aim is to get contradiction from (\ref{AplbA}), (\ref{Bprime})
and (\ref{Aprime}).

Let us use the symbol
\begin{equation}\label{Kdef}
K=\max_{b\in
B^{'}}|A^{'}+bA^{'}|\phantom{111111}\mbox{so}\phantom{111111}K<
\frac{2^{219}\sqrt{2}\ln^{16}(e|A|)}{C^{61}}|B|^{61\delta}|A|.
\end{equation}
Now we use Lemma \ref{XYGlemma} to establish that

\begin{multline*}
E_{+}(A^{'},bA^{'})=|\{(a_1,a_2,a_3,a_4)\in A^{'}\times A^{'}\times
A^{'}\times A^{'}:
a_1+a_2 b = a_3+a_4 b \}|\geqslant\\
\geqslant\frac{|A^{'}|^{4}}{|A^{'}+b
A^{'}|}\geqslant\frac{|A^{'}|^4}{K}.
\end{multline*}
Summing over all $b\in B^{'}$ we obviously obtain
$$
|\{(a_1,a_2,a_3,a_4,b)\in A^{'}\times A^{'}\times A^{'}\times
A^{'}\times B^{'}: a_1+a_2 b = a_3+a_4 b
\}|\geqslant\frac{|A^{'}|^{4}|B^{'}|}{K}.
$$
There are some elements $\widetilde{a}_2, \widetilde{a}_3\in A^{'}$
such that
$$
|\{(a_1,a_4,b)\in A^{'}\times A^{'}\times B^{'}:
a_1-\widetilde{a}_3= (a_4-\widetilde{a}_2)b
\}|\geqslant\frac{|A^{'}|^{2}|B^{'}|}{K}.
$$
Let $A^{'}_1=A^{'}-\widetilde{a}_3, A^{'}_2=A^{'}-\widetilde{a}_2$
be translates of $A^{'}$ by $\widetilde{a}_3$ and $\widetilde{a}_2$
respectively. Then
$$
|\{(a_1,a_2,b)\in A^{'}_1\times A^{'}_2\times B^{'}: a_1= a_2 b
\}|\geqslant\frac{|A^{'}|^{2}|B^{'}|}{K}.
$$
There is some $a_{*}\in A_2^{'}$ such that
$$
|\{(a_1,b)\in A^{'}_1\times B^{'}: a_1= a_{*} b
\}|\geqslant\frac{|A^{'}||B^{'}|}{K}.
$$
Thus, we have a subset $B^{'}_1\subset (A^{'}_1\cap a_{*}B^{'})$ of
cardinality
$$|B^{'}_1|\geqslant\frac{|A^{'}||B^{'}|}{K}.$$
In original notations $B^{'}_1$ lies in the intersection of
$\frac{a_{*}}{b_0}B_2$ and some translate of $A_1^{(b_0)}$; besides
by the bounds (\ref{Bprime}), (\ref{Aprime}) and (\ref{Kdef})
\begin{equation}\label{Boneprime}
|B^{'}_1|>\frac{C^{69}}{2^{250}\ln^{18}(e|A|)}|B|^{1-69\delta}.
\end{equation}

We consider three cases.

1) Case 1. Suppose that $Q[B_1^{'}]\neq \Fp$. It is clear that
$1+Q[B_1^{'}]\neq Q[B_1^{'}]$ since otherwise $Q[B_1^{'}]=\Fp$. The
latter mean that there are elements $a,b,c,d\in B_1^{'}$ with
$1+\frac{a-b}{c-d}\notin Q[B_1^{'}]$. Now we recall that $B_1^{'}$
is a subset of $\frac{a_{*}}{b_0}B_2$ so we can regard $a,b,c,d$ as
elements of $B_2$. Observe, that for an arbitrary subset
$B_{1}^{''}\subset B_1^{'}, |B_{1}^{''}|\geqslant 0.98 |B_{1}^{'}|$
we have $1+\frac{a-b}{c-d}\notin Q[B_{1}^{''}]$ since
$Q[B_{1}^{''}]\subset Q[B_{1}^{'}]$. Therefore, by Lemma
\ref{QXYproplemma}, for these elements $a,b,c,d\in B_2$ we have
\begin{equation}\label{tmplab8}
(0.98)^2|B_1^{'}|^2\leqslant |B_1^{''}|^2=
\left|B_1^{''}+\left(B_1^{''}+\frac{a-b}{c-d}B_1^{''}\right)\right|\leqslant
\left|B_1^{''}+B_1^{''}+\frac{a-b}{c-d}B_1^{''}\right|.
\end{equation}

We now use Lemma \ref{ChenCovering}. Let us first show that for any
$b_1\in B_2$ we can cover 99\% of the elements of the set
$b_1B_1^{'}$ (a subset of the translation of $b_1A_1^{(b_0)}$) or
$-b_1B_1^{'}$ by at most
$\frac{2^{109}\ln(100)\ln^8(e|A|)}{C^{28}}|B|^{28\delta}$ additive
translates of the set $b_0A_1^{(b_0)}$. Indeed $b_0A_1^{(b_1)}\cap
A_1^{(b_0)}$ is a subset of $b_0A_1^{(b_0)}$, and by Lemma
\ref{ChenCovering} and Lemma \ref{RuszaTriangle}) we can cover 99\%
of the elements of either $b_1B_1^{'}$ or $-b_1B_1^{'}$ by at most
$$\frac{\ln(100)}{|b_0A_1^{(b_1)}\cap
A_1^{(b_0)}|}\min\left\{|b_0A_1^{(b_1)}\cap
A_1^{(b_0)}+b_1B_1^{'}|,|b_0A_1^{(b_1)}\cap
A_1^{(b_0)}-b_1B_1^{'}|\right\}\leqslant$$
$$\leqslant\frac{\ln(100)}{|A_1^{(b_1)}\cap
A_1^{(b_0)}|}\min\left\{|b_0A_1^{(b_1)}\cap
A_1^{(b_0)}+b_1A_1^{(b_0)}|,|b_0A_1^{(b_1)}\cap
A_1^{(b_0)}-b_1A_1^{(b_0)}|\right\}\leqslant $$
$$\leqslant\frac{\ln(100)|A_1^{(b_1)}\cap A_1^{(b_0)}+b_1A_2^{(b_0)}\cap
A_2^{(b_1)}||A_1^{(b_0)}+b_0A_2^{(b_0)}\cap
A_2^{(b_1)}|}{|A_1^{(b_1)}\cap A_1^{(b_0)}||b_0b_1A_2^{(b_1)}\cap
A_2^{(b_0)}|}\leqslant$$
$$\leqslant\frac{\ln(100)|A_1^{(b_1)}+b_1A_2^{(b_1)}||A_1^{(b_0)}+b_0A_2^{(b_0)}|}
{|A_1^{(b_1)}\cap A_1^{(b_0)}||A_2^{(b_1)}\cap
A_2^{(b_0)}|}\leqslant
\frac{2^{105}\ln(100)\ln^8(e|A|)}{C^{28}}|B|^{28\delta}$$additive
translates of $b_0A_1^{(b_1)}\cap A_1^{(b_0)}$ and whence of
$b_0A_1^{(b_0)}$. In the last estimate we have used
(\ref{tmplabQb}), (\ref{tmplab5}) and (\ref{tmplab7}).

This altogether enables us to choose $B^{''}_1$ as a subset
containing at least 98\% of the elements from $B_1^{'}$ such that
$(a-b)B^{''}_1$ gets covered by at most
$\frac{2^{210}\ln^2(100)\ln^{16}(e|A|)}{C^{56}}|B|^{56\delta}$
translates of $b_0A_1^{(b_0)}+b_0A_1^{(b_0)}$. Similarly, we can
find a subset $\widetilde{A}_1^{(b_0)}$ containing at least 98\% of
the elements of $A_1^{(b_0)}$ such that
$(c-d)\widetilde{A}_1^{(b_0)}$ gets covered by at most
$\frac{2^{210}\ln^2(100)\ln^{16}(e|A|)}{C^{56}}|B|^{56\delta}$
translates of $b_0A_1^{(b_0)}+b_0A_1^{(b_0)}$. Now we apply Lemma
\ref{RuszaGeneral} to (\ref{tmplab8}) as follows
\begin{multline}\label{tmplab9}
    \left|B_1^{''}+B_1^{''}+\frac{a-b}{c-d}B_1^{''}\right|\leqslant
    \frac{|\widetilde{A}_1^{(b_0)}+B_1^{''}+B_1^{''}||
    \widetilde{A}_1^{(b_0)}+\frac{a-b}{c-d}B_1^{''}|}{|\widetilde{A}_1^{(b_0)}|}\leqslant\\
    \leqslant \frac{2^4\sqrt{2}|B|^{\delta}}{ C|A|}|A_1^{(b_0)}+A_1^{(b_0)}+A_1^{(b_0)}|
    |\widetilde{A}_1^{(b_0)}+\frac{a-b}{c-d}B_1^{''}|\leqslant\\
    \leqslant\frac{2^{87}\ln^{6}(e|A|)}{C^{25}}|B|^{25\delta}|\widetilde{A}_1^{(b_0)}+\frac{a-b}{c-d}B_1^{''}|
\end{multline}

The covering arguments above implies that
$$|\widetilde{A}_1^{(b_0)}+\frac{a-b}{c-d}B_1^{''}|\leqslant
\frac{2^{420}\ln^{4}(100)\ln^{32}(e|A|)}{C^{112}}|B|^{112\delta}
|A_1^{(b_0)}+A_1^{(b_0)}+A_1^{(b_0)}+A_1^{(b_0)}|\leqslant
$$
$$\leqslant\frac{2^{530}\ln^{4}(100)\ln^{40}(e|A|)}{C^{144}}|B|^{144\delta}|A|.$$
Comparing to (\ref{Boneprime}) and using the condition
$\frac{|B|}{|A|}\geqslant\frac{1}{4}$, for large $p$ we deduce
$$\frac{(0.98)^2C^{138}}{2^{500}\ln^{36}(e|A|)}|B|^{2-138\delta}<
\frac{2^{613}\ln^{4}(100)\ln^{46}(e|A|)}{C^{169}}|B|^{169\delta}|A|\Leftrightarrow$$
\begin{equation}\label{TmplabB}
\Leftrightarrow\frac{|B|^{2-307\delta}}{|A|\ln^{82}(e|A|)}<\frac{2^{1113}\ln^{4}(100)}
{(0,98)^2C^{307}}\Rightarrow
|B|^{1-308\delta}<\frac{2^{1115}\ln^{4}(100)} {(0,98)^2C^{307}}.
\end{equation}
Now we define
$C=\frac{2^{\frac{1115}{307}}\ln^{\frac{4}{307}}(100)}{(0.98)^{\frac{2}{307}}}$
and from (\ref{TmplabB}) deduce the inequality
$$|B|<|B|^{308\delta}$$
which is false when $\delta\leqslant\frac{1}{308}$. This finishes
proof of the Theorem \ref{main1} in case 1.

2) Case $2$. Suppose that $|B_1^{'}|>\sqrt{p}$. It is clear that
$Q[B_1^{'}]=\Fp$ since for an arbitrary $\xi\in\Fp$ the equality
$|B_1^{'}+\xi B_1^{'}|=|B_1^{'}|^2$ is impossible (simply because
$|B_1^{'}|^2>p$). Let us take arbitrary elements $\xi\in\Fp^{*}$,
$s\in\Fp$, an arbitrary subset $|B_1^{''}|\geqslant 0.96|B_1^{'}|$
and denote
$$f_{\xi}(s):=|\{(b_1,b_2)\in B_1^{'}\times B_1^{'}:b_1+\xi b_2=s\}|$$
$$f_{\xi}^{'}(s):=|\{(b_1,b_2)\in B_1^{''}\times B_1^{''}:b_1+\xi b_2=s\}| $$

It is obvious that
$$\sum_{s\in
\Fp}(f_{\xi}(s))^{2}=|\{(b_1,b_2,b_3,b_4)\in B_1^{'}\times
B_1^{'}\times B_1^{'}\times B_1^{'}: b_1+\xi b_2=b_3+\xi b_4\}|$$
$$=|B_1^{'}|^2+|\{(b_1,b_2,b_3,b_4)\in B_1^{'}\times
B_1^{'}\times B_1^{'}\times B_1^{'}:b_1\neq b_3, b_1+\xi b_2=b_3+\xi
b_4\}|$$ and
$$
\sum_{s\in\Fp}f_{\xi}(s)=|B_1^{'}|^2
$$
$$\sum_{s\in\Fp}f_{\xi}^{'}(s)=|B_1^{''}|^2.$$

Let us observe that for every $b_1,b_2, b_3, b_4\in B_1^{'}$ such
that $b_1\neq b_3$, there is at most one $\eta\in\Fp^{*}$ satisfying
the equality $b_1+\eta b_2=b_3+\eta b_4$. Therefore,
$$\sum_{\xi\in\Fp^{*}}\sum_{s\in\Fp}(f_{\xi}(s))^{2}\leqslant
|B_1^{'}|^2(p-1)+|B_1^{'}|^{4}.$$ From the last inequality it
directly follows that there is an element $\xi\in\Fp^{*}$ such that
$$
\sum_{s\in\Fp}(f_{\xi}^{'}(s))^{2}\leqslant
\sum_{s\in\Fp}(f_{\xi}(s))^{2}\leqslant
|B_1^{'}|^2+\frac{|B_1^{'}|^{4}}{p-1}.
$$
Note that this $\xi$ is independent on $B_1^{''}$. According to
Cauchy-Schwartz,
$$
\left(\sum_{s\in\Fp}f_{\xi}^{'}(s)\right)^{2}\leqslant |B_1^{''}+\xi
B_1^{''}|\sum_{s\in\Fp}(f_{\xi}^{'}(s))^{2}.
$$
Now we see that
\begin{equation}\label{tmplab20}
|B_1^{''}+\xi B_1^{''}|\geqslant
\frac{|B_1^{''}|^4(p-1)}{|B_1^{'}|^2(p-1)+|B_1^{'}|^4}\geqslant
\frac{(0.96)^4|B_1^{'}|^4(p-1)}{|B_1^{'}|^2(p-1)+|B_1^{'}|^4}\geqslant
(0.96)^4\frac{p-1}{2}.
\end{equation}

Reminding that $Q[B^{'}_1]=\Fp$, we can find elements $a,b,c,d\in
B_1^{'}$, such that $\xi=\frac{a-b}{c-d}$ (again, we can regard them
as elements of $B_2$). Using similar covering arguments as in proof
of the case 1 we can deduce that we can choose $B^{''}_1$ as a
subset containing at least 96\% of the elements from $B_1^{'}$ such
that $(a-b)B^{''}_1+(c-d)B^{''}_1$ gets covered by at most
$\frac{2^{420}\ln^{4}(100)\ln^{32}(e|A|)}{C^{112}}|B|^{112\delta}$
translates of
$b_0A_1^{(b_0)}+b_0A_1^{(b_0)}+b_0A_1^{(b_0)}+b_0A_1^{(b_0)}$. Now
we see that
$$
\left|B_1^{''}+\frac{a-b}{c-d}B_1^{''}\right|\leqslant
\frac{2^{420}\ln^{4}(100)\ln^{32}(e|A|)}{C^{112}}|B|^{112\delta}
|A_1^{(b_0)}+A_1^{(b_0)}+A_1^{(b_0)}+A_1^{(b_0)}|\leqslant$$
$$\leqslant\frac{2^{530}\ln^{4}(100)\ln^{40}(e|A|)}{C^{144}}|B|^{144\delta}|A|.
$$
Again, comparing to (\ref{tmplab20}) and using the condition
$\frac{|B|}{|A|}\geqslant \frac{1}{4}$, we deduce
$$(0.96)^4\frac{p}{4}\leqslant(0.96)^4\frac{p-1}{2}<
\frac{2^{530}\ln^{4}(100)\ln^{40}(e|A|)}{C^{144}}|B|^{144\delta}|A|\Rightarrow$$
\begin{equation}\label{TmplabB2}
\Rightarrow\frac{p}{4}<\frac{2^{530}\ln^{4}(100)}{C^{144}(0.96)^4}p^{145\beta\delta+\alpha}
\end{equation}
Now we define
$C=\frac{2^{\frac{265}{72}}\ln^{\frac{1}{36}}(100)}{(0.96)^{\frac{1}{36}}}$
and from (\ref{TmplabB2}) deduce the inequality
$$p<p^{145\beta\delta+\alpha}$$
which is false when $\delta\leqslant\frac{1-\alpha}{145\beta}$. This
concludes proof of the Theorem in case 2.

3) Case $3$. Suppose that $Q[B_1^{'}]=\Fp$ and $|B_1^{''}|\leqslant
\sqrt{p}$. Repeating arguments from the proof of case 2 for an
arbitrary subset $B_1^{''}\subset B_1^{'}, |B_1^{''}|\geqslant
0.96|B_1^{'}|$ we finding elements $a,b,c,d\in B_2$ independent on
the subset $B_1^{''}$ with
$$\left|B_1^{''}+\frac{a-b}{c-d}B_1^{''}\right|\geqslant(0.96)^4\frac{|B_1^{'}|^2}{2}.$$

Using similar covering arguments as in proof of the case 1 we can
deduce that we can choose $B^{''}_1$ as a subset containing at least
96\% of the elements from $B_1^{'}$ such that
$(a-b)B^{''}_1+(c-d)B^{''}_1$ gets covered by at most
$\frac{2^{420}\ln^{4}(100)\ln^{32}(e|A|)}{C^{112}}|B|^{112\delta}$
translates of
$b_0A_1^{(b_0)}+b_0A_1^{(b_0)}+b_0A_1^{(b_0)}+b_0A_1^{(b_0)}$. Now
we see that
$$
\left|B_1^{''}+\frac{a-b}{c-d}B_1^{''}\right|\leqslant
\frac{2^{420}\ln^{4}(100)\ln^{32}(e|A|)}{C^{112}}|B|^{112\delta}
|A_1^{(b_0)}+A_1^{(b_0)}+A_1^{(b_0)}+A_1^{(b_0)}|\leqslant$$
$$\leqslant\frac{2^{530}\ln^{4}(100)\ln^{40}(e|A|)}{C^{144}}|B|^{144\delta}|A|.
$$
Comparing to (\ref{Boneprime}) and using the condition
$\frac{|B|}{|A|}\geqslant\frac{1}{4}$, we deduce
$$\frac{(0.96)^4C^{138}}{2^{500}\ln^{36}(e|A|)}|B|^{2-138\delta}<
\frac{2^{530}\ln^{4}(100)\ln^{40}(e|A|)}{C^{144}}|B|^{144\delta}|A|\Leftrightarrow$$
\begin{equation}\label{TmplabB1}
\Leftrightarrow\frac{|B|^{2-282\delta}}{|A|\ln^{76}(e|A|)}<\frac{2^{1030}\ln^{4}(100)}
{(0,96)^4C^{282}}\Rightarrow
|B|^{1-283\delta}<\frac{2^{1032}\ln^{4}(100)} {(0,96)^4C^{282}}.
\end{equation}
Now we define
$C=\frac{2^{\frac{516}{141}}\ln^{\frac{2}{141}}(100)}{(0.96)^{\frac{2}{141}}}$
and from (\ref{TmplabB1}) deduce the inequality
$$|B|<|B|^{283\delta}$$
which is false when $\delta\leqslant\frac{1}{283}$. Note that in all
the cases the meaning assigned for the constant $C$ is strictly less
than $15$. The Theorem \ref{main1} is proved. $\blacksquare$
\bigskip

\section{Proof of the Theorem \ref{mainLargeB}.}\label{sec3}

As in the proof of the Proposition \ref{main1} we assume contrary,
i.e.
$$\sum_{b\in B}E_{+}(A, bA)>C|B|^{1-\delta}|A|^3$$
for some $C>0$, $\delta>0$. Following arguments in the beginning of
the proof of the Proposition \ref{main1}, we finding $A^{'}\subset
A$ and $B^{'}\subset\Fp^{*}, 1\in B^{'}$ (which is in fact a subset
of a multiplicative shift of $B$) such that
\begin{equation}\label{tmplabel3}
|A^{'}+bA^{'}|<\frac{2^{219}\sqrt{2}\ln^{16}
(e|A|)}{C^{61}}|B|^{61\delta}|A|=K\phantom{1}\textmd{for}\phantom{1}\textmd{all}\phantom{1}b\in
B^{'}
\end{equation}
\begin{equation}\label{tmplabel4}
|B^{'}|>\frac{C^7}{2^{26}\ln^2(e|A|)}|B|^{1-7\delta}
\end{equation}
\begin{equation}\label{tmplabel5}
|A^{'}|>\frac{C}{2^{4}\sqrt{2}}|B|^{-\delta}|A|.
\end{equation}
Using Lemma \ref{XYGlemma} we obtain
\begin{multline*}
|\{(a_1,a_2,a_3,a_4)\in A^{'}\times A^{'}\times A^{'} \times
A^{'}:a_1+ba_2=a_3+ba_4\}|>\\
>\frac{|A^{'}|^4}{K}\phantom{1}\mbox{for all $b\in B^{'}$}.
\end{multline*}
Summing up by all $b\in B^{'}$ one gets
\begin{multline*}
|\{(a_1,a_2,a_3,a_4,b)\in A^{'}\times A^{'}\times A^{'} \times
A^{'}\times B^{'}:a_1+ba_2=a_3+ba_4\}|>\\
>\frac{|A^{'}|^4|B^{'}|}{K}\phantom{1}\mbox{for all $b\in
B^{'}$}.
\end{multline*}
Now we can fix elements $a_3^{0}, a_2^{0}\in A^{'}$ such that
\begin{equation}\label{tmplabel21}
|\{(a_1,a_4,b)\in A^{'}\times A^{'}\times
B^{'}:a_1-a_3^{0}=b(a_4-a_2^{0})\}|
>\frac{|A^{'}|^2|B^{'}|}{K}.
\end{equation}
We denote
$$f(s)=|\{(a,b)\in A^{'}\times B^{'}:b(a-a_2^{0})=s\}|,$$
$$g(s)=\left\{%
\begin{array}{ll}
    1, & \hbox{if $s\in A^{'}-a_3^{0}$;} \\
    0, & \hbox{otherwise.} \\
\end{array}%
\right.$$ Clearly,
\begin{equation}\label{tmplabel1}
|\{(a_1,a_4,b)\in A^{'}\times A^{'}\times
B^{'}:a_1-a_3^{0}=b(a_4-a_2^{0})\}|=\sum_{s\in\Fp}f(s)g(s),
\end{equation}
\begin{equation}\label{tmplabel2}
\sum_{s\in\Fp}f^2(s)=E_{\times}(A^{'}-a_2^{0},B^{'}).
\end{equation}
Now, by Cauchy-Schwartz,
$$\left(\sum_{s\in\Fp}f(s)g(s)\right)^2\leqslant \sum_{s\in\Fp}f^{2}(s)\sum_{s\in\Fp}g^{2}(s)$$
and, by (\ref{tmplabel1}) and (\ref{tmplabel2}), one can deduce
$$E_{\times}(A^{'}-a_2^{0},B^{'})>\frac{|A^{'}|^3|B^{'}|}{K^{2}}.$$

Consider two cases.

Case 1. Assume that $|A^{'}||B^{'}|\leqslant p$. Applying Theorem
\ref{garaev} one obtains
$$\frac{K^4}{|A^{'}|}>|A^{'}-A^{'}|^2\cdot\frac{|A^{'}|^2|B^{'}|^2}{E_{\times}(A^{'}-a_2^{0},B^{'})}\geqslant C_1
\frac{|A^{'}|^3|B^{'}|^{\frac{1}{9}}}{\log_2(|B^{'}|)}.$$ Using
(\ref{tmplabel3}), (\ref{tmplabel4}) and (\ref{tmplabel5}) we deduce
$$\frac{C_1C^{\frac{43}{9}}|B|^{\frac{1}{9}-\frac{43}{9}\delta}|A|^{4}}
{2^{\frac{188}{9}}\ln^{\frac{2}{9}}(e|A|)\log_2(|B|)}<\frac{2^{878}\ln^{64}(e|A|)}{C^{244}}|B|^{244\delta}|A|^4\Rightarrow
$$
\begin{equation}\label{tmplabel6}
|B|^{\frac{1}{9}}<\frac{2^{\frac{8090}{9}}\ln^{\frac{578}{9}}(e|A|)\log_2(|B|)}{C_1C^{\frac{2239}{9}}}|B|^{\frac{2239}{9}\delta}.
\end{equation}
Defining $C=\frac{2^{\frac{8090}{2239}}}{C_1^{\frac{9}{2239}}}$, we
observe that for sufficiently large $p$ from (\ref{tmplabel6})
follows the inequality
$$|B|^{\frac{1}{9}}<|B|^{\frac{2240}{9}\delta}.$$
which gives a contradiction when $\delta=\frac{1}{2240}.$ This
completes proof of the Theorem \ref{mainLargeB}  in this case.

Case 2. Assume that $|A^{'}||B^{'}|> p$. Again, applying Theorem
\ref{garaev} we obtain
$$\frac{K^4}{|A^{'}|}>|A^{'}-A^{'}|^2\cdot\frac{|A^{'}|^2|B^{'}|^2}{E_{\times}(A^{'}-a_2^{0},B^{'})}\geqslant C_1
\frac{|A^{'}|^{\frac{26}{9}}p^{\frac{1}{9}}}{\log_2p}.$$ Using
(\ref{tmplabel3}) and (\ref{tmplabel5}) we deduce
$$\frac{C_1C^{\frac{35}{9}}|A|^{\frac{35}{9}}p^{\frac{1}{9}}}
{2^{\frac{35}{2}}|B|^{\frac{35}{9}}\log_2p}<\frac{2^{878}\ln^{64}(e|A|)}{C^{244}}|B|^{244\delta}|A|^4\Rightarrow
$$
\begin{equation}\label{tmplabel7}
\Rightarrow\frac{2^{\frac{1791}{2}}\ln^{64}(e|A|)\log_2p}{C^{\frac{2231}{9}}C_1}
|A|^{\frac{1}{9}}|B|^{\frac{2231}{9}\delta}>p^{\frac{1}{9}}.
\end{equation}
Defining $C=\frac{2^{\frac{19119}{4462}}}{C_1^{\frac{9}{2231}}}$, we
observe that for sufficiently large $p$ from (\ref{tmplabel7})
follows the inequality
$$p^{\frac{1}{9}}<|B|^{\frac{2232}{9}\delta}|A|^{\frac{1}{9}}.$$
which gives a contradiction when $\delta=\frac{1-\alpha}{2232}.$
Theorem \ref{mainLargeB} is proved. $\blacksquare$


\begin{thebibliography}{}

\bibitem{B} J. Bourgain, \emph{Multilinear exponential sums in prime fields
under optimal entropy condition on the sources}, Geometric and
Functional Analysis, vol. 18, N 5, 2009 , pp. 1477 -- 1502.

\bibitem{BoGl} J. Bourgain, A. Glibichuk, \emph{Exponential sum estimate over subgroup in an
arbitrary finite field}, accepted for publication in Journal de
Analyze Math\'{e}matiques.

\bibitem{BDFKK} J. Bourgain, S. J. Dilworth, K. Ford, S. Konyagin, D.
Kutzarova, \emph{Explicit constructions of RIP matrices}, Proc. 43rd
ACM Symposium of the Theory of Computing (STOC), pp. 637---644
(2011).

\bibitem{RH} H. Helfgott, M. Rudnev, \emph{An explicit incidence theorem in
$\Fp$}, preprint, arXiv:1001.1980v2.

\bibitem{R1} I. Z. Ruzsa, \emph{An application of graph theory to additive number
theory}, Scientia, Ser. A, 3 (1989), 97 -- 109.

\bibitem{R2} I. Z. Ruzsa, \emph{Sums of finite sets}, Number theory (New York,
1991 -- 1995), 281 -- 293, Springer, New York, 1996.

\bibitem{BG} J. Bourgain and M. Z. Garaev, \emph{On a variant of sum-product estimates and
explicit exponential sums bounds in prime fields}, Mathematical
proceedings of the Cambridge Philosophical Society, vol. 146 (2009),
part 1, pp. 1 -- 21.

\bibitem{TV} T. Tao and V. Vu, \emph{Additive combinatorics},
Cambridge University Press, Cambridge, 2006.

\bibitem{BKT} J.~Bourgain, N.~Katz, T.~Tao, \emph{A sum-product estimate in
finite fields and their applications}, Geom and Funct. Anal., {\bf
14} (2004), 27--57.

\bibitem{BK} J.~Bourgain, S.~Konyagin, \emph{Estimates for the number
of sums and products and for exponential sums over subgroups in
fields of prime order}, C.R. Acad. Sci. Paris, Ser. I, {\bf 337}
(2003), 75--80.

\bibitem{Shen} Chun-Yen Shen, \emph{Quantitative sum product estimates on different sets},
Electron. J. Combin., 15 (2008), no. 1.

\bibitem{Gar} M. Z. Garaev, \emph{Sums and products of sets and estimates of rational
trigonometric sums in fields of prime order}, Russian Mathematical
Surveys, 2010, vol. 65, no. 4, pp. 599--658


\end{thebibliography}
\end{document}